\documentclass{article}
\usepackage{texdraw}

\newtheorem{thm}{Theorem}[section]

\newtheorem{prop}[thm]{Proposition}

\newtheorem{lem}[thm]{Lemma}

\newtheorem{rem}[thm]{Remark}

\newtheorem{ex}{Example}[section]

\newcommand{\be}{\begin{equation}}
\newcommand{\ee}{\end{equation}}
\newcommand{\ben}{\begin{enumerate}}
\newcommand{\een}{\end{enumerate}}
\newcommand{\pa}{{\partial}}

\newcommand{\R}{{\rm R}}
\newcommand{\e}{{\epsilon}}

\newcommand{\pxi}{ {\pa \over \pa x^i}}
\newcommand{\pxj}{ {\pa \over \pa x^j}}

\font\BBb=msbm10 at 12pt
\newcommand{\Bbb}[1]{\mbox{\BBb #1}}

\newcommand{\qed}{\hspace*{\fill}Q.E.D.}  
\title{\Large Finsler Metrics with ${\bf K}=0$ and ${\bf S}=0$}
\author{Zhongmin Shen}

\date{May,  2001; Revised in September, 2001}

\begin{document}
\maketitle

\begin{abstract}
{In Finsler geometry, there are infinitely many models of constant curvature. The Funk metrics, the Hilbert-Klein metrics and the Bryant metrics are projectively flat with non-zero constant curvature. 
A recent  example constructed by the author is  projectively flat 
 with zero curvature. 
In this paper, we introduce a technique to construct non-projectively flat  Finsler metrics with zero curvature in each dimension. The technique can be used to construct many non-projectively flat Finsler metrics of constant curvature.
}
\end{abstract}

\section{Introduction}

One of the fundamental problems in Riemann-Finsler geometry is  to study and characterize Finsler metrics of constant flag curvature. 
The flag curvature in Finsler geometry is an analogue of the sectional curvature in Riemannian geometry. 
Riemannian metrics of constant  sectional curvature were classified  by E. Cartan a long time ago.
There are only three local Riemannian metrics of constant sectional curvature, up to a scaling. 
However,  the local metric structure of a Finsler metric with constant flag curvature is much more complicated.
Mathematicians have discovered several important  Finsler metrics of constant flag curvature (e.g.,
the Funk metrics (${\bf K}=-1/4$), the Hilbert-Klein metrics (${\bf K}=-1$)\cite{BaChSh}\cite{Sh1} and the Bryant metrics (${\bf K}=1$)\cite{Br1}\cite{Br2}\cite{Br3}). All of them are locally projectively flat.

Recently, 
the author and D. Bao discover a family of non-projectively flat 
Finsler metrics of constant curvature ${\bf K}=1$ on $\Bbb S^3$ using the Hopf-fibration structure \cite{BaSh}.
We know that 
Minkowski spaces are  Finsler spaces with ${\bf K}=0$.
In \cite{Sh2}, we construct a positively complete  projectively flat spray with vanishing Riemann curvature (R-flat). 
Then we construct a positively complete  locally projectively flat
Finsler metric with ${\bf K}=0$. It is defined on the unit ball $\Bbb B^n\subset \Bbb R^n$ by
\be
F({\bf y}) := { \Big ( \sqrt{ |{\bf y}|^2 - (|{\bf x}|^2|{\bf y}|^2 
-\langle {\bf x}, {\bf y} \rangle^2) } + \langle {\bf x}, {\bf y} \rangle \Big )^2 
\over ( 1-|{\bf x}|^2)^2 \sqrt{|{\bf y}|^2 - (|{\bf x}|^2|{\bf y}|^2 
-\langle {\bf x}, {\bf y} \rangle^2)} }, \ \ \ \ \ \ {\bf y}\in T_{\bf x}\Bbb B^n=\R^n,\label{fromFunk}
\ee 
where $|\cdot |$ and $\langle \; , \; \rangle$ denote the standard Euclidean 
norm and inner product. According to \cite{AZ}, any positively complete Finsler metric with 
${\bf K}=0$ must be locally Minkowskian if the  first and second Cartan torsions are bounded (see also \cite{Sh1}). Thus the Finsler metric 
in (\ref{fromFunk}) does not have bounded Cartan torsions.

 In this paper, we are going to construct a non-projectively flat Finsler metrics with ${\bf K}=0$ in each dimension. Further, these Finsler metrics have bounded
first and second Cartan torsions. Thus, they are not positively complete.

\begin{thm}\label{thmK=0}
Let $n\geq 2$ and
\[ \Omega:=\Big \{ p=(x, y, \bar{p}) \in \R^2 \times \R^{n-2} 
\ \Big | \ x^2+y^2 < 1 \Big \} .\]
Define
\be
 F({\bf y}) := {\sqrt{ \Big ( -y u + x v \Big )^2 
+ |{\bf y}|^2 \Big ( 1- x^2 -y^2 \Big ) } - \Big ( - yu + x v \Big ) 
\over 1- x^2 - y^2 },\label{RandersExample}
\ee
where ${\bf y} = (u, v,\bar{\bf y} ) \in T_p \Omega = \R^n$ and $p= (x, y, \bar{p})\in \Omega$.
$F$ is a Finsler metric on $\Omega$ satisfying 
\[ {\bf K} =0, \ \ \ \ \ \ {\bf S}=0.\]
\end{thm}

Here the quantity ${\bf S}$ denotes the S-curvature which is a non-Riemannian quantity. It is known that every  Berwald metric satisfies  ${\bf S}=0$ (\cite{Sh1}\cite{Sh4}). 
Since every Berwald metric with ${\bf K}=0$ must be
locally Minkowskian (see \cite{AIM}\cite{BaChSh}), the Finsler metric in (\ref{RandersExample}) is not Berwaldian.

The construction of the Finsler metric in (\ref{RandersExample}) is motivated by the shortest time problem (see Section 3 below).

The Finsler metric in (\ref{RandersExample}) is in the following form
$
 F = \alpha + \beta, $
where $\alpha $ is a Riemannian metric and $\beta $ is a 1-form with $\|\beta\|_{\alpha}(x):=\sup_{y\in T_xM} \beta(y)/\alpha(y) <1$ for all $x\in M$. 
This type of Finsler metrics were first studied by G. Randers in 1941 \cite{Ra} from the standard point of general relativity (see also  \cite{AIM}). Therefore they are called  {\it Randers metrics}.

The Randers metric in (\ref{RandersExample}) is not locally projectively 
flat, hence not locally Minkowskian. Theorem \ref{thmK=0} is inconsistent  with 
the main result in 
\cite{SSAY}, where they claim that every Randers metric with ${\bf K}=0$ must be locally Minkowskian. See also \cite{YaSh}\cite{Ma} for related discussion. 
 Nevertheless, if a Randers metric with ${\bf K}=0$ is positively complete, then it must be
locally Minkowskian.

\begin{thm}\label{thm1}
Let $F=\alpha+\beta$ be a positively complete Randers metric on a manifold $M$.
Then ${\bf K}=0$ if and only if it is locally Minkowskian.
In this case, $\alpha$ is flat and $\beta$ is parallel with respect to $\alpha$.
\end{thm}

Here  a Finsler metric is said to be positively complete if every geodesic defined on $(a, b)$ can be extended to a geodesic defined on $(a, \infty)$. Example (\ref{RandersExample}) shows that
the positive completeness can not be dropped. 
In \cite{Sh6}, we show that every locally projectively flat Randers metric 
with ${\bf K}=0$ must be locally Minkowskian. In 
this case, the positive completeness is not required.

\section{Preliminaries}

In this section, we recall some basic definitions in Riemann-Finsler geometry \cite{BaChSh}\cite{Sh1}.

A Finsler metric on a manifold $M$ is a function $F: TM \to [0, \infty)$ 
with the following properties:
\ben
\item[(a)] $F(\lambda {\bf y}) = \lambda F({\bf y}) $, $\lambda >0$;
\item[(b)] For any non-zero vector ${\bf y}\in T_xM$,
the induced bilinear form  $g_{\bf y}$ on $T_xM$ is an inner product, where
\[ g_{\bf y}({\bf u}, {\bf v}):= {1\over 2} {\pa^2\over\pa s\pa t}\Big [ F^2({\bf y}+s {\bf u}+t {\bf v}) \Big ] |_{s=t=0}, \ \ \ \ \ \ {\bf u}, {\bf v}\in T_xM.
\]
\een

\bigskip
Riemannian metrics are special Finsler metrics. Traditionally, a Riemannian metric is denoted by $ a_{ij}(x) dx^i \otimes dx^j$. It is a family of inner products on tangent spaces. Let
 $\alpha({\bf y}):= \sqrt{g_{ij}(x)y^iy^j}$,
${\bf y}= y^i\pxi|_x\in T_xM$. $\alpha$ is  a family of Euclidean norms on tangent spaces.   Throughout this 
paper, we also denote   a Riemannian metric by  $\alpha=\sqrt{a_{ij}(x)y^iy^j}$.

Let $\alpha=\sqrt{a_{ij}(x)y^iy^j}$ be a Riemannian metric  and $\beta =b_i(x)y^i$ a 1-form on a manifold $M$. Define
\[ F = \alpha +\beta.\]
Then $F$ satisfies (a). If we assume that 
\[ \|\beta\|_{\alpha}(x)
:= \sup_{y\in T_xM}{\beta ({\bf y})\over\alpha({\bf y}) }=\sqrt{ a^{ij}(x) b_i(x)b_j(x)} <1 , \ \ \ \ \ \ \ x\in M,\]
then $F$ satisfies (b). By definition, $F$ is a Finsler metric.
We call $F$ a Randers metric.

\bigskip

Let $F$ be a Finsler metric. For a non-zero vector ${\bf y}\in T_pM \setminus\{0\}$,
define 
\[ {\bf C}_{\bf y} ({\bf u}, {\bf v}, {\bf w} ) 
:= {1\over 4 } {d^3 \over dt^3} \Big [ F^2 \Big ( {\bf y}+ s {\bf u}+ t {\bf v}+ r {\bf w} \Big ) \Big ]|_{s=t=r=0}\]
\[ \tilde{\bf C}_{\bf y}  (  {\bf u}, {\bf v}, {\bf w}, {\bf z} ) 
:= {1\over 4 } {d^4 \over dt^4} \Big [ F^2 \Big ( {\bf y}+ s {\bf u}+ t {\bf v}+ r {\bf w}+ h{\bf z} \Big ) \Big ]|_{s=t=r=h=0}.\]
By the homogeneity of $F$, we have
\[ {\bf C}_{\bf y} ({\bf u}, {\bf v}, {\bf y} ) =0, \ \ \ \ \
\tilde{\bf C}_{\bf y}  (  {\bf u}, {\bf v}, {\bf w}, {\bf y}  ) 
= - {\bf C}_{\bf y} ({\bf u}, {\bf v}, {\bf w} ).\]
We call ${\bf C}$ and $\tilde{\bf C}$ the first and second Cartan torsion respectively. We usually call ${\bf C}$ the Cartan torsion if the second one does not appear.
It is obvious that $F$ is Riemannian if and only if ${\bf C}=0$,
and ${\bf C}=0$ if and only if $\tilde{\bf C}=0$.

The essential bounds of ${\bf C}$ and $\tilde{\bf C}$ at $x\in M$ are defined by
\begin{eqnarray*}
\|{\bf C}\|_x : & = &  \sup_{ {\bf u} \in T_xM} F({\bf y}){ |{\bf C}_{\bf y} ({\bf u}, {\bf u}, {\bf u} \Big ) |\over  [g_{\bf y}({\bf u}, {\bf u} ) ]^{3/2} }\\
\|\tilde{\bf C}\|_x:& = &  \sup_{ {\bf u} \in T_xM}F^2({\bf y}) { |\tilde{\bf C}_{\bf y} ({\bf u}, {\bf u}, {\bf u}, {\bf u} \Big ) |\over  [g_{\bf y}({\bf u}, {\bf u} ) ]^2 },
\end{eqnarray*}
where the supremum is taken over all non-zero vectors ${\bf y}, {\bf u}\in T_xM$ with $g_{\bf y}({\bf y}, {\bf u}) =0$.

\bigskip
For a Finsler metric $F$ on an $n$-dimensional manifold $M$, the (Busemann-Hausdorff) volume form $dV_F = \sigma_F(x) dx^1 \cdots dx^n$ is defined by
\be
\sigma_F(x) := {{\rm Vol} (\Bbb B^n(1)) \over {\rm Vol} \Big \{ (y^i)\in \R^n \ \Big | \ F \Big ( y^i \pxi|_x \Big ) < 1 \Big \} } .\label{dV}
\ee

In general, the local scalar function 
$\sigma_F(x)$ can not be expressed in terms of elementary functions, even $F$ is locally expressed by elementary functions. However, for Randers metrics, the volume form is expressed by a very simple formula, since each indicatrix 
$S_xM := \{ {\bf y}\in T_xM \ | \ F({\bf y}) = 1 \}$ is a shifted Euclidean sphere
in $(T_xM, \alpha_x)$.
More precisely,
for a Randers metric 
$ F = \alpha + \beta$, where $\alpha =\sqrt{a_{ij}(x)y^iy^j}$ 
and $\beta(y)=b_i(x)y^i$, its volume form $dV_F$ is given by
\be
dV_F = \Big ( 1 - \|\beta \|^2_{\alpha}(x) \Big )^{n+1\over 2} dV_{\alpha}.\label{dV_F}
\ee
 See \cite{Sh1} for more details.

\bigskip
Now we introduce the notion of S-curvature for Finsler metrics
(see \cite{Sh1}).
Let $F$ be a Finsler metric on a manifold $M$.
Let \[ g_{ij}(x,y):= g_{\bf y} \Big (\pxi|_x, \pxj|_x\Big )= {1\over 2}[F^2]_{y^iy^j}(x,y), \ \ \ \ \ {\bf y}= y^i \pxi|_x.\]
Express the volume form $dV_F$ by
\[ dV_F = \sigma_F(x) dx^1 \cdots dx^n,\]
where $\sigma_F(x)$ is defined in (\ref{dV}).
Set
\[ \mu({\bf y}):= \ln \Big [ {\sqrt{\det (g_{ij}(x, y)) }\over \sigma_F(x) }
\Big ].\]
The quantity $\mu$ is a scalar function on $TM\setminus\{0\}$.
$\mu$ is a non-Riemannian quantity. 

\bigskip

Locally minimizing constant speed curves (geodesics) are characterized by
\[ {d^2 x^i\over dt^2} + 2 G^i(x, {dx\over dt} ) =0,\]
where $G^i(x,y)$ are given by
\be
G^i := {1\over 4} g^{il} \Big \{ 2 {\pa g_{jl}\over \pa x^k} - {\pa g_{jk}
\over \pa x^l} \Big \} y^j y^k.
\ee
$G^i$ are called the geodesic coefficients in a local coordinate system.
If $F$ is Riemannian, then $G^i(x, y)= {1\over 2} \Gamma^i_{jk}(x)y^jy^k$ 
are quadratic in $(y^i)$ at every point $x\in M$.
A Finsler metric is called a {\it Berwald metric} if the geodesic coefficients
have this property. There are many non-Riemannian Berwald metrics.
The classification of Berwald metrics is done by Z.I. Szabo \cite{Sz}.

\bigskip

For a vector ${\bf y}\in T_xM$
Let $c(t), -\e < t <\e $, denote the geodesic with $c(0)=x$ and $\dot{c}(0)={\bf y}$. Define
\[ {\bf S}({\bf y}):= {d \over dt} \Big [ \mu \Big (\dot{c}(t)\Big ) \Big ]|_{t=0}.\]
We call ${\bf S}$ the S-curvature. This quantity was introduced in \cite{Sh4}\cite{Sh1}.

Let $G^i(x, y)$ denote the geodesic coefficients of $F$ in the same local coordinate system. The S-curvature is defined by
\be
{\bf S}({\bf y}) := {\pa G^i\over \pa y^i}(x,y) - y^i {\pa \over \pa x^i} \Big [ \ln \sigma_F (x)\Big ],\label{Slocal}
\ee
where ${\bf y}=y^i\pxi|_x\in T_xM$.
This quantity was first introduced in \cite{Sh4} for a volume comparison theorem.
 It is proved that ${\bf S}=0$ if $F$ is a Berwald metric \cite{Sh4}.
There are many non-Berwald metrics satisfying ${\bf S}=0$.

\bigskip

Now, we  recall the definition of 
Riemann curvature. Let $F$ be a Finsler metric on an $n$-manifold
and $G^i$ denote the geodesic coefficients of $F$.
For a vector ${\bf y}=y^i\pxi|_x\in T_xM$, define
${\bf R}_{\bf y}= R^i_{\ k}(x, y) dx^k \otimes \pxi|_x: T_xM \to T_xM$ by
\be
 R^i_{\ k} := 2 {\pa G^i\over \pa x^k} 
- y^j {\pa^2 G^i\over \pa x^j\pa y^k} + 2 G^j {\pa^2 G^i\over \pa y^j \pa y^k} 
- {\pa G^i\over \pa y^j} {\pa G^j\over \pa y^k}.
\label{Rik}
\ee
The Ricci curvature is defined by
\[
{\bf Ric} ({\bf y}):= R^i_{\ i} (x, y).
\]
In dimension two, let $x:=x^1, y:=x^2,  u: =y^1, v:= y^2$.
We can express the Ricci curvature by
\begin{eqnarray}
{\bf Ric}({\bf y}) & = & 2 \Big \{ {\pa G^1 \over \pa x} + {\pa G^2 \over \pa y} + {\pa G^1 \over \pa u} {\pa G^2 \over \pa v} 
- {\pa G^1 \over \pa v} {\pa G^2 \over \pa u} \Big \}\nonumber\\
&&  - S^2 -\Big ( u {\pa \over x} + v {\pa \over \pa y} - 2 G^1 {\pa \over \pa u} - 2 G^2 {\pa \over \pa v} \Big ) (S), \label{RicRic}
\end{eqnarray}
where $S = {\pa G^1 \over \pa u} + {\pa G^2 \over \pa v}$.

\bigskip

For a two-dimensional tangent plane $P\subset T_xM$ and a non-zero vector 
${\bf y}\in P$, define
\be
{\bf K}(P, {\bf y}):= { g_{\bf y} \Big ( {\bf R}_{\bf y}({\bf u}),\; {\bf u} \Big )\over g_{\bf y}({\bf y}, {\bf y}) g_{\bf y}({\bf u}, {\bf u})
- g_{\bf y}({\bf y}, {\bf u}) g_{\bf y} ({\bf y}, {\bf u})},\label{flag}
\ee
where $P= {\rm span}\{ {\bf y}, {\bf u}\}$.  ${\bf K}$ is called the flag curvature. 
Usually, ${\bf K}(P, {\bf y})$ depends on the direction ${\bf y}\in P$.
If $F$ is Riemannian, then ${\bf K}(P, {\bf y})$ is independent of ${\bf y}\in P$. 
The flag curvature is an analogue of the sectional curvature in Riemannian geometry. 
It is obvious that ${\bf K}=0$ if and only if ${\bf R}=0$.

The following theorem is due to Akbar-Zadeh \cite{AZ}.
\bigskip
\begin{thm}\label{thmAZ} {\rm (\cite{AZ})}
Let $(M, F)$ be a positively complete Finsler manifold with ${\bf K}=0$.
Assume that 
\[ \sup_{x\in M} \|{\bf C}\|_x <\infty, \ \ \ \ \ \sup_{x\in M} \|\tilde{\bf C}\|_x < \infty.\]
Then $F$ is locally Minkowskian.
\end{thm}

See also Theorem 10.3.7 in \cite{Sh1} for a proof.

\bigskip

Let $F$ be a Finsler metric on a manifold $M$ and $G^i (x, y)$ denote the geodesic coefficients of $F$. 
Let $\tilde{F}$ be another metric and $\tilde{G}^i(x,y)$ denote the geodesic coefficients of $\tilde{F}$. To find the relationship between
the Riemann curvature $R^i_{\ k}(x,y)dx^k \otimes \pxi $ of $F$ and the Riemann curvature $\tilde{R}^i_{\ k}(x,y) dx^k \otimes \pxi$ of $\tilde{F}$, we introduce 
\[ H^i(x, y): = G^i(x, y) - \tilde{G}^i(x,y).\]
Define
\[ H^i_{\ | k}: = {\pa H^i\over \pa x^k} + H^j {\pa^2 \tilde{G}^i\over \pa y^j \pa y^k}
- {\pa H^i\over \pa y^j} {\pa \tilde{G}^j \over \pa y^k}.\]
We have the following useful formula
\be
 R^i_{\ k} = \tilde{R}^i_{\ k} + 2 H^i_{\; |k} - y^j (H^i_{\; |j})_{y^k} + 2 H^j (H^i)_{y^jy^k} - (H^i)_{y^j}(H^j)_{y^k}.\label{RRH}
\ee
The proof is straight forward, so is omitted. See \cite{Sh1}.

\section{Shortest Time Problem}
Let $(M, F )$ be a Finsler space.
Suppose that an object on $(M, F)$ is pushed by an internal force ${\bf u}$  
with constant length,  $F({\bf u}) =1$. 
Due to the friction, the object moves on $M$ at a constant speed, but it can change direction freely.  
Without  external force acting on the object, any path of shortest time is a shortest 
path of $F$.

Now given an external force ${\bf v}$ pushing the object.  We assume  that $F(-{\bf v}) < 1$, otherwise, the object can not move forward in the direction  $-{\bf v}$.
Due to friction, the speed of the object is proportional  to
the length of the combined force 
\[ {\bf t}:= {\bf v} + {\bf u}.\]
For simplicity, we may assume that 
the speed of the object  is equal to $F({\bf t})$. 
That is, ${\bf t}:={\bf v}+{\bf u}$ is the velocity vector once the direction of the internal force is chosen.

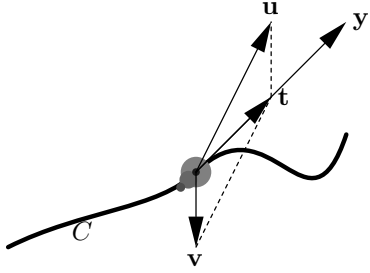
\begin{center}
\begin{texdraw}
\drawdim cm
\linewd  0.06

\move(-2 0)
\clvec(-1 0.5)(0 0.5)(0.5 1)
\clvec(1.5 2)(2 0)(2.5 1.5)

\move(-1 0.2)
\textref h:C v:C \htext{$C$}

\linewd 0.02

\move(0.5 1)
\fcir f:0.5 r:0.2
\move(0.4 0.9)
\fcir f:0.4 r:0.12
\move(0.3 0.8)
\fcir f:0.3 r:0.06

\move(0.5 1)
\fcir f:0.1 r:0.05
\arrowheadtype t:F \avec(1.5 2)
\move(0.5 1)
\arrowheadtype t:F \avec(0.5 0)
\move(0.5 1) 
\arrowheadtype t:F \avec(1.5 3)
\move(1.5 3.1)
\textref h:C v:B \htext{${\bf u}$}
\move(0.5 -0.1)
\textref h:C v:T \htext{${\bf v}$}
\move(1.6 2)
\textref h:L v:C \htext{${\bf t}$}
\move(0.5 1)
\arrowheadtype t:F \avec(2.5 3)
\move(2.6 3)
\textref h:L v:C \htext{${\bf y}$}

\lpatt(0.05 0.06)
\move(1.5 2)
\lvec(1.5 3)
\move(0.5 0)
\lvec(1.5 2)

\move(0 3.5)
\fcir f:1 r:0.0
\move(0 -0.5)
\fcir f:1 r:0.0
\end{texdraw}
\end{center}

Express the velocity vector in the following  form
\[ {\bf t} = F({\bf t}) {\bf y},\]
where ${\bf y}$ is a unit vector with respect to $F$. Since $F({\bf u})=1$,
$F({\bf t})$ is determined by
\be
 F \Big ( F ({\bf t}) {\bf y} - {\bf v} \Big ) =F({\bf u}) = 1.\label{eqv}
\ee

Now we are going to find the Finsler metric $\tilde{F}$
such that the $\tilde{F}$-length of any  curve is equal to 
the time for which the object travels along it.

\begin{lem} Let $(M, F)$ and ${\bf v}$ be a vector field on $M$ with $F(-{\bf v}) < 1$.
 Define $\tilde{F}: TM \to [0, \infty)$ by
\be
F \Big ( {{\bf y}\over \tilde{F}({\bf y}) } - {\bf v} \Big ) =1,
\ \ \ \ \ \ {\bf y}\in T_xM\setminus\{0\}.
\label{tF}
\ee
 For any curve $C$ in $M$, the $\tilde{F}$-length of $C$ is equal to
the time for which the object travels along it.
\end{lem}
{\it Proof}: Take an arbitrary curve $C$ from $p$ to $q$,  
and a coordinate map $c: [0, T] \to C$ such that 
$c(0)=p$,  $c(T)=q$ and the velocity vector $\dot{c}(t)$ is equal to the 
combined force  at $c(t)$. 
Express $\dot{c}(t) =F(\dot{c}(t)) {\bf y}(t)$,
where ${\bf y}(t)$ is tangent to $C$ at $c(t)$ with $F({\bf y}(t))=1$.
By (\ref{eqv}),
\be 
F \Big ( F (\dot{c}(t)) {\bf y}(t) - {\bf v}(t) \Big ) =1,\label{lam1}
\ee
where ${\bf v}(t):= {\bf v}_{c(t)}$. 
Consider the following equation
\[ F\Big ( \lambda {\bf y}(t) - {\bf v}(t) \Big ) =1, \ \ \ \ \ \lambda >0.\]
Since $F(-{\bf v}(t)) < 1$, the above equation has a unique solution, that is $F(\dot{c}(t))$. By the definition of $\tilde{F}$, we obtain
\be
\tilde{F}({\bf y}(t)) = {1\over F ( \dot{c}(t) ) }. \label{lam2}
\ee
From (\ref{lam1}) and (\ref{lam2}), we obtain
\[ \tilde{F}(\dot{c}(t)) 
= F (\dot{c}(t)) \tilde{F}({\bf y}(t)) = 1.\]
This implies 
\[  T = \int_0^T \tilde{F}(\dot{c}(t)) dt.\]
Thus $\tilde{F}$ is the desired Finsler metric.
\qed
\bigskip

The relation between $F$ and $\tilde{F}$ is actually very simple. 
For a point $x\in M$, the indicatrix of $F$ at $x$ is related to 
that of $\tilde{F}$ at $x$, 
\[ \Big \{ {\bf y}\in T_xM \ | \ \tilde{F}({\bf y}) = 1 \Big \}
= \Big \{ {\bf y}\in T_xM \ | \ F({\bf y} ) = 1 \Big \} + {\bf v}_x.\]
This leads to the following

\begin{prop}\label{lemdvdv} Let 
$(M, F)$ be a Finsler manifold and ${\bf v}$ be a vector field on $M$ with
$F(-{\bf v}) < 1$. Let  $\tilde{F}: TM \to [0, \infty)$ denote the
Finsler metric defined by (\ref{tF}).
The volume form of $F$ equals that of $\tilde{F}$, 
\[ dV_F = dV_{\tilde{F}}.\]
\end{prop}
{\it Proof}: 
Let $(x^i)$ be a local coordinate system at $x\in M$
and 
\begin{eqnarray*}
{\cal U}_F : & = & \Big \{ (y^i)\in \R^n \ | \ F\Big ( y^i\pxi|_x \Big ) < 1 \Big \},\\
{\cal U}_{\tilde{F}} : & = & \Big \{ (y^i)\in \R^n \ | \ F\Big ( y^i\pxi|_x \Big ) < 1 \Big \}.
\end{eqnarray*}
From (\ref{tF}), we have
\[ {\cal U}_{\tilde{F}} = {\cal U}_F + (v^i),\]
where ${\bf v} = v^i\pxi|_x$. 
This gives that 
\[ dV_F = {{\rm Vol}(\Bbb B^n(1))\over {\rm Vol} ({\cal U}_F )}\; dx^1\cdots dx^n = {{\rm Vol} (\Bbb B^n(1)) \over {\rm Vol} ({\cal U}_{\tilde{F}})}\; dx^1\cdots dx^n = dV_{\tilde{F}}.\]
\qed 

\bigskip

\begin{rem}
{\rm Consider the shortest time problem on a Riemannian manifold 
$(M, \alpha)$ with an external force field ${\bf v}$. 
Denote by $\langle \; , \; \rangle_{\alpha}$  the family of  inner products on tangent spaces, which are  determined by $\alpha$,
\[ \alpha({\bf y}) =\sqrt{\langle{\bf y}, {\bf y}\rangle_{\alpha}}, \ \ \ \ \ {\bf y}\in T_pM.\]
Solving the following equation for $\tilde{F}({\bf y})$,
\[\alpha \Big ( {{\bf y}\over \tilde{F}({\bf y}) } - {\bf v} \Big ) = 1,\]
we obtain
\be
 \tilde{F} ({\bf y}) =\tilde{\alpha}({\bf y}) + \tilde{\beta}({\bf y}), \ \ \ \ \ \ {\bf y}\in T_pM\label{Randersmetric}
\ee
where
\[
\tilde{\alpha}({\bf y}):=
{ \sqrt{\langle {\bf v}, {\bf y}\rangle_{\alpha}^2 + \alpha({\bf y})^2 
(1-\alpha({\bf v})^2) }\over 1 -\alpha({\bf v})^2 },
\ \ \ \ \ \tilde{\beta}({\bf y}) := - {\langle {\bf v}, {\bf y}\rangle_{\alpha}
\over 1-\alpha({\bf v})^2 }. 
\]
By Proposition \ref{lemdvdv}, we have
\[ dV_{\tilde{F}} = dV_{\alpha}.\]
}
\end{rem}

The Finsler metric $\tilde{F}$ in (\ref{Randersmetric}) is called a {\it Randers metric} in Finsler geometry.

\bigskip

\begin{ex}{\rm
Let $\alpha$
denote the standard Euclidean metric on the unit ball $\Bbb B^n$ and
 ${\bf v}$ denote  the radial  vector field on
$\Bbb B^n$, which is given by
\[ {\bf v}_p = - (x^i), \ \ \ \ \ \ \ \ p=(x^i)\in \Bbb B^n.\]
The Randers metric  associated with $(\alpha, {\bf v})$ as defined in (\ref{Randersmetric}) is given by
\[ \tilde{F} ({\bf y}) = { \sqrt{ \langle {\bf v}, {\bf y} \rangle^2 + |{\bf y}|^2 (1-|{\bf v}|^2 ) } + \langle {\bf v}, {\bf y} \rangle \over 1- |{\bf v}|^2}, 
\ \  \ \ {\bf y}\in T_p\Bbb B^n,\]
where $|\cdot |$ and $\langle\; , \; \rangle$ denote the standard
Euclidean norm and inner product.
This is just the Funk metric on $\Bbb B^n$. 
Geodesics of $\tilde{F}$ are straight lines. Moreover, the flag curvature is a
negative constant, ${\bf K}= -1/4$. 
If an object moves away from the center, it takes infinite time to reach the boundary. However, it  takes finite time to reach the center along any shortest path.  Thus $\tilde{F}$ is positively complete, but not negatively complete.
}
\end{ex}

\begin{ex}\label{exBaSh} {\rm 
(\cite{BaSh}) Let $\Bbb S^3$ be the standard unit sphere in $\Bbb R^4$.
Let $\alpha$ denote the standard Riemannian metric on $\Bbb S^3$ and 
 ${\bf v}: = \e {\bf w}$, where $ |\e | < 1$ and ${\bf w}$ is  a left-invariant unit vector field on $\Bbb S^3$. 
The Randers metric associated with $(\alpha, {\bf v})$ as defined in (\ref{Randersmetric})
is given by
\be
\tilde{F} ({\bf y} )
= { \sqrt{ \e^2 \langle {\bf w}, {\bf y} \rangle^2_{\alpha} + (1-\e^2) \alpha ({\bf y})^2    } - \e \langle {\bf w}, {\bf y} \rangle_{\alpha} \over 1- \e^2 },
\ \ \ \ \ \ \ {\bf y}\in T_p\BBb S^3.
\ee
We have shown that ${\bf K}=1$ for any $\e$ with $|\e | < 1$. 
}
\end{ex}

Other examples will be discussed in Section 7 below  and \cite{Sh5}.
\bigskip

\section{Proof of Theorem \ref{thm1}}

According to Theorem \ref{thmAZ}, 
to prove Theorem \ref{thm1}, it suffices to prove that the first and second Cartan torsions have uniform upper bounds.

\bigskip

Consider a Randers metric $ F =\alpha +\beta$,
where $ \alpha =\sqrt{a_{ij}y^iy^j}$ and $\beta=b_i y^i$ 
with $ \|\beta\|_{\alpha}  =\sqrt{a^{ij}b_i b_j} < 1$.
It is proved that the first Cartan torsion of $F$ satisfies the bound
\be
 \|{\bf C}\|_x \leq {3\over \sqrt{2}} \sqrt{ 1-\sqrt{ 1-\|\beta\|^2(x)}} < {3\over \sqrt{2}}, \ \ \ \ \  \ x\in M.\label{CBound}
\ee
This is verified by B. Lackey in dimension two \cite{BaChSh}, and can be extended to higher dimensions with a
simple argument \cite{Sh3}. 

Now we are going to find an upper bound on the second Cartan torsion for Randers metrics.
First, we consider a special two-dimensional case.

Let  $p\in M$
and $\kappa:=\|\beta\|_{\alpha}(p) < 1$. There is an  orthonormal basis $\{ {\bf e}_1, {\bf e}_2 \}$ for $T_pM$ such that 
\[ F ({\bf y})=\sqrt{u^2+v^2} + \kappa u, \ \  \ \ \ \ {\bf y}= u{\bf e}_1 + v {\bf e}_2 \in T_pM.\]

\begin{lem}\label{lemsecondCartan}
Let $\kappa$ be a constant with $ 0 \leq \kappa < 1$ and
\be
F := \sqrt{u^2+v^2} + \kappa u.
\ee
The second Cartan torsion of $F$ satisfies the following bound
\[ \|\tilde{\bf C}\| \leq {27\over 2} \kappa <  {27\over 2} .\]
\end{lem}
{\it Proof}: Let ${\bf y}=  r (\cos \theta, \sin \theta) $, where $ r >0$, and ${\bf y}^{\bot}$ denote 
the vector perpendicular to ${\bf y}$ with respect to $g_{\bf y}$,
\[ g_{\bf y} ({\bf y}, {\bf y}^{\bot} ) =0, \ \ \ \ \ \
g_{\bf y} ({\bf y}^{\bot}, {\bf y}^{\bot} ) = F^2({\bf y}).\]
We have
\[ {\bf y}^{\bot}= {r\over \sqrt{ 1+\kappa \cos \theta} } ( -\sin \theta, \ \kappa + \cos\theta ).\]
By a direct computation, we obtain 
\[ \tilde{\bf C}_{\bf y} ({\bf y}^{\bot},{\bf y}^{\bot},{\bf y}^{\bot},{\bf y}^{\bot}) = F^2 ({\bf y})\Big \{ 
6 \kappa { \kappa +\cos \theta \over 1+ \kappa \cos \theta } - {15\over 2} \kappa \cos \theta \Big \}.\]
This gives
\begin{eqnarray*} \|\tilde{\bf C}\|
& = &  \max_{ 0 \leq \theta \leq 2\pi}  \Big |6 \kappa { \kappa +\cos \theta \over 1+ \kappa \cos \theta } - {15\over 2} \kappa \cos \theta \Big | \\
&  \leq & \max_{ 0 \leq \theta \leq 2\pi}\Big \{ 6 \kappa\Big | { \kappa +\cos \theta \over  1 + \kappa \cos \theta } \Big | + {15\over 2} \kappa \Big |\cos \theta \Big |\Big \}\\
& \leq & 6 \kappa + {15\over 2}\kappa = {27\over 2} \kappa .
\end{eqnarray*}
This proves the lemma. \qed

\bigskip
We can extend Lemma \ref{lemsecondCartan} to higher dimensions. 
\begin{lem}\label{lem3.3}
Let $F= \alpha +\beta$ be a Randers metric on an $n$-manifold $M$. 
The second Cartan torsion of $F$ satisfies
\[ \|\tilde{\bf C}\|\leq {27\over 2} \|\beta\|_{\alpha} < {27\over 2}.\]
\end{lem}
{\it Proof}:
At a point $p\in M$, there are two vectors ${\bf y}, {\bf u}\in T_pM$
with 
\[ F({\bf y}) =1, \ \ \ \ g_{\bf y}({\bf y}, {\bf u}) =0, \ \ \ \ g_{\bf y}({\bf u}, {\bf u}) =1\]
 such that
\[ \|\tilde{\bf C}\|_p = |\tilde{\bf C}_{\bf y} ({\bf u}, {\bf u}, {\bf u}, {\bf u})|
.\]
Let $V:={\rm span} \{ {\bf y}, {\bf u} \}$ and 
\[ \kappa := \sup_{{\bf v}\in V} {\beta({\bf v})\over \alpha({\bf v})}.\]
By Lemma \ref{lemsecondCartan}, 
\[ |\tilde{\bf C}_{\bf y} ({\bf u}, {\bf u}, {\bf u}, {\bf u})| \leq {27\over 2} \kappa.\]
Note that 
\[ \kappa =\sup_{{\bf v}\in V} {\beta({\bf v})\over \alpha({\bf v})}\leq \sup_{{\bf v}\in T_pM}{\beta({\bf v})\over \alpha({\bf v})} = : \|\beta\|_{\alpha} (p).\]
We obtain 
\[ \|\tilde{\bf C}\|_p \leq {27\over 2} \|\beta \|_{\alpha}(p).\]
\qed

\bigskip
\noindent
{\it Proof of Theorem \ref{thm1}}: By (\ref{CBound}),
we know that
$\|{\bf C}\| < 3/\sqrt{2}$. By (\ref{lem3.3}), we know that
$\|\tilde{\bf C}\| < 13.5$.  Then Theorem \ref{thm1} follows from Theorem \ref{thmAZ}.
\qed

\section{Randers Metrics with ${\bf S}=0$}

In this section, we are going to  find a sufficient and necessary condition on $\alpha$ and $\beta$ for ${\bf S}=0$.
In particular, we will show that if
$\beta$ is a Killing form of constant length, then 
${\bf S}=0$.

Let $F =\alpha +\beta$ be a Randers metric on a manifold $M$,
where 
\[\alpha (y)=\sqrt{a_{ij}(x) y^iy^j}, \ \ \ \ \ \beta(y)= b_i(x) y^i\]
with 
$\|\beta\|_x :=\sup_{y \in T_xM} \beta(y)/\alpha(y) < 1$. 

In a standard local coordinate system $(x^i, y^i)$ in $TM$,
define $b_{i|j}$ by
\[ b_{i|j} \theta^j := db_i -b_j \theta_i^{\ j},\]
where 
$\theta^i :=dx^i$ and 
$\theta_i^{\ j} :=\tilde{\Gamma}^j_{ik} dx^k$ denote the Levi-Civita  connection forms of $\alpha$.
Let 
\[ r_{ij}:= {1\over 2} \Big ( b_{i|j} + b_{j|i} \Big ), \ \ \ \ 
s_{ij} := {1\over 2} \Big ( b_{i|j}-b_{j|i} \Big ),\]
\[ s^i_{\ j} =a^{ip}s_{\ pj} \ \ \ \ s_j:= b_is^i_{\ j} .\]
The geodesic  coefficients $G^i$ of $F$ are related to the geodesic coefficients $\tilde{G}^i$ of $\alpha$ by 
\be
G^i = \tilde{G}^i + P y^i + Q^i ,\label{Gi}
\ee
where 
\begin{eqnarray*}
P: & =& {1\over 2F} \Big \{ r_{ij} y^iy^j -2\alpha  s_i y^i \Big \}\\
Q^i: & = & \alpha s^i_{\ j}y^j.
\end{eqnarray*}
See \cite{AIM}.
Observe that  
\[ {\pa Q^i \over \pa  y^i} = \alpha^{-1} y_i s^i_{\ j}y^j+ \alpha s^i_{\ i} =\alpha^{-1}s_{ij}y^iy^j + \alpha a^{ij}s_{ij}=
0,\]
where $y_i :=a_{ij}y^j$.
Thus
\[ {\pa G^i\over \pa y^i} = {\pa \tilde{G}^i\over \pa y^i}  + (n+1) P.\]
Put 
\[ dV_F :=\sigma_F(x) dx^1 \cdots dx^n, \ \ \ \ dV_{\alpha} = \sigma_{\alpha}(x) dx^1\cdots dx^n.\]
 According to (\ref{Gi}), 
we have 
\[ \sigma_F = \Big ( 1 -\|\beta\|^2_{\alpha} \Big ) ^{n+1\over 2} \sigma_{\alpha}.\]
Note that 
\[ d \Big [ \ln \sigma_{\alpha} \Big ] ={\pa \tilde{G}^i\over \pa y^i}.\]
By (\ref{Slocal}), we obtain a formula for ${\bf S}$,
\be
{\bf S} = (n+1) \Big \{ P - d \Big [ \ln \sqrt{ 1-\|\beta\|^2_{\alpha} } \Big ] \Big \}. \label{SRandersS}
\ee

We have the following
\begin{prop}\label{propS=c}
Let $F=\alpha+\beta$ be a Randers metric on an $n$-manifold $M$,
where $\alpha=\sqrt{a_{ij}(x)y^iy^j}$ and $\beta =b_i(x) y^i$.
Then
\be
{\bf S}=0\label{ScF}
\ee 
 if and only if 
\be
r_{ij} +b_is_j + b_j s_i = 0.\label{eqS=0}
\ee
\end{prop}
{\it Proof}: 
For the sake of simplicity, we
choose an orthonormal basis for $T_xM$ such that $a_{ij}=\delta_{ij}$.
Let
\[ \rho := \ln\sqrt{1-\|\beta\|_{\alpha}^2 }.\]
and
$d\rho =\rho_i dx^i$, i.e.,
\be
\rho_i = - { b_j b_{j|i}\over 1-\|\beta\|_{\alpha}^2 }.\label{rhoi}
\ee
By (\ref{SRandersS}), 
${\bf S}=0$ if and only if  
\be
r_{ij} y^i y^j - 2\alpha s_i y^i = 2 (\alpha+\beta) \rho_i y^i.
\label{eqS=0*}
\ee
(\ref{eqS=0*}) is equivalent to 
the following 
equations
\begin{eqnarray}
&& r_{ij} = b_j \rho_i + b_i \rho_j  \label{eqrij}\\
&& -s_i = \rho_i \label{eqspl}
\end{eqnarray}

First we assume that 
${\bf S}=0$. Then (\ref{eqrij}) and (\ref{eqspl}) hold.
Plugging (\ref{eqspl}) into (\ref{eqrij}) gives (\ref{eqS=0}).

Now we assume that (\ref{eqS=0}) holds. 
Note that $s_jb_j = b_i s_{ij} b_j =0$. 
Contracting (\ref{eqS=0}) with $b_j$ yields
\be
 b_j r_{ij} =  - \|\beta\|^2_{\alpha} s_i , \label{eqsss}
\ee
that is,
\[ b_j b_{i|j} + b_j b_{j|i} =- \|\beta \|^2 (b_j b_{j|i} - b_j b_{i|j} ).\]
We obtain
\be
 b_j b_{i|j} =  - {1+\|\beta\|^2_{\alpha} 
\over 1-\|\beta\|^2_{\alpha} } \; b_j b_{j|i}.\label{bijij}
\ee
It follows from  (\ref{rhoi}) and (\ref{bijij}) that
\be
 s_i =  {1\over 2} \Big ( b_j b_{j|i}
+ {1+\|\beta\|^2_{\alpha} 
\over 1-\|\beta\|^2_{\alpha} } \; b_j b_{j|i} \Big ) =  - \rho_i.
\label{s_i}
\ee
Thus
\begin{eqnarray*}
r_{ij} y^iy^j - 2\alpha s_i y^i & = &
 - 2 \beta s_i y^i - 2\alpha s_i y^i\\
& = & - 2 (\alpha+\beta ) s_i y^i \\
& = &  2 F \rho_i y^i.
\end{eqnarray*}
We obtain
\[
{\bf S}= (n+1) \Big \{\rho_i y^i - \rho_i y^i \Big \}
= 0.
\]
This gives (\ref{ScF}).\qed

\section{Randers Metrics with ${\bf K}=0$ and ${\bf S}=0$}

In this section, we are going to compute the Riemann curvature
of a Randers metric satisfying ${\bf S}=0$.

Let $F=\alpha+\beta$ be a Randers metric on a manifold $M$. 
Denote by $b_{i|j}$, $b_{i|j|k}$, etc the coefficients of the covariant derivatives of $\beta$ with respect to $\alpha$. Set
$r_{ij}:= (b_{i|j}+b_{j|i})/2$, $s_{ij}:=(b_{i|j}-b_{j|i})/2$, $s^i_{\ j}=a^{ik}s_{kj}$ 
and $s_j := b_is^i_{\ j}$. Further,  we set $s_0:=s_py^p, s^i_{\ 0}:=s^i_{\ p} y^p $, $s_{0|j}=s_{p|j}y^p$  and $s_{0|0}=s_{p|q}y^py^q$, etc.
Denote the Riemann curvature of $\alpha$ by 
$\tilde{R}^i_{\ k} dx^k \otimes \pxi$. We have the following

\begin{thm}\label{thmK=0S=0}
Let $F = \alpha +\beta$ be a Randers metric satisfying ${\bf S}=0$. Then 
${\bf K}=0$   if and only if the following two equations hold,
\begin{eqnarray}
\bar{R}^i_{\ k} & = & -\Big ( s_{0|0} \delta^i_k -s_{0|k} y^i \Big )
- \Big ( s_{k|0}-s_{0|k} \Big )y^i \nonumber\\
& & - s_0\Big ( s_0\delta^i_k- s_k y^i\Big ) +\Big (\alpha^2 s^i_{\ j}s^j_{\ k}-  s^i_{\ j}s^j_{\ 0} y_k \Big ) -3 s_{k0} s^i_{\ 0}, \label{K=0C}\\
0 & = & s_js^j_{\ 0} \Big ( \alpha^2 \delta^i_k - y_k y^i \Big )+
\alpha^2 \Big ( s_j s^j_{\ 0} \delta^i_k - s_j s^j_{\ k} y^i \Big )\nonumber\\
&& + \alpha^2 \Big ( s^i_{\ k|0}-s^i_{\ 0|k} \Big )
- \Big ( \alpha^2 s^i_{\ 0|k}-s^i_{\ 0|0} y_k  \Big ). \label{K=0B} 
\end{eqnarray}
\end{thm}
{\it Proof}:
By assumption ${\bf S}=0$, we obtain from (\ref{SRandersS}) and (\ref{eqspl}), we obtain
\[ P= - s_0.\] Thus
$ G^i = \tilde{G}^i + H^i,$
where 
\[ H^i = -s_0 y^i + \alpha s^i_{\ 0}.\]
By a direct computation, we obtain
\begin{eqnarray*}
H^i_{\; |k} & = & -s_{0|k} y^i +\alpha s^i_{\ 0|k}  \\
\Big ( H^i_{\; |j} \Big )_{y^k} & = & -s_{k|j} y^i -s_{0|j} \delta^i_k  + \alpha^{-1} y_k s^i_{\ 0|j}  + \alpha s^i_{\ k|j}\\
y^j \Big ( H^i_{\; |j} \Big )_{y^k}
& = & -s_{k|0} y^i -s_{0|0} \delta_k^i   + \alpha^{-1} y_k s^i_{\ 0|0} +\alpha s^i_{\ k|0} \\
\Big (H^i\Big )_{y^j} & = & 
-s_j y^i -s_0\delta^i_j 
+ \alpha^{-1} y_j s^i_{\ p}y^p +\alpha s^i_{\ j}\\
\Big ( H^i\Big )_{y^jy^k} & = & -s_j \delta^i_k -s_k \delta^i_j 
+\alpha^{-3}(a_{jk}\alpha^2 - y_j y_k) s^i_{\ 0} 
 + \alpha^{-1}( y_j s^i_{\ k} 
+  y_k s^i_{\ j}),
\end{eqnarray*}
$y_k :=a_{jk}y^j$.
Plugging them into (\ref{RRH}), we obtain
\begin{eqnarray}
R^i_{\ k} &=& \bar{R}^i_{\ k}+\Big ( s_{0|0} \delta^i_k -s_{0|k} y^i \Big )
+ \Big ( s_{k|0}-s_{0|k} \Big )y^i
\nonumber\\
&&+ s_0\Big ( s_0\delta^i_k- s_k y^i\Big ) -\Big (\alpha^2 s^i_{\ j}s^j_{\ k}-  s^i_{\ j}s^j_{\ 0} y_k \Big ) +3 s_{k0} s^i_{\ 0}\nonumber\\
&& - \Big \{ s_js^j_{\ 0} \Big ( \alpha^2 \delta^i_k - y_k y^i \Big )
+\alpha^2 \Big ( s_j s^j_{\ 0} \delta^i_k - s_j s^j_{\ k} y^i \Big )\nonumber\\
&& +\alpha^2 \Big ( s^i_{\ k|0}-s^i_{\ 0|k} \Big )
- \Big ( \alpha^2 s^i_{\ 0|k}-s^i_{\ 0|0} y_k  \Big ) \Big \}\alpha^{-1}
\label{localRik}
\end{eqnarray}

According to (\ref{localRik}), the coefficients of the Riemann curvature
of $F$ are in the following form
\[ R^i_{\ k} = A + B\alpha^{-1} ,\]
where $A$ and $B$ are polynomials of $y^i$ at each point $x\in M$. 
Thus $R^i_{\ k}=0$ if and only if $A=0$ and $B=0$.
This proves Theorem \ref{thmK=0S=0}.
\qed

\bigskip

Taking the trace of $R^i_{\; k}$ in (\ref{localRik}), we obtain 
\begin{eqnarray}
{\bf Ric} & = & \overline{\bf Ric} 
+ (n-1) \Big \{s_{0|0} + s_0s_0 \Big \} + 2 s_{k0}s^k_{\ 0} -\alpha^2 s^k_{\ j} s^j_{\ k} \nonumber \\
&& + 2\Big \{  s^k_{\ 0|k} - (n-1)  s_js^j_{\ 0}   \Big \}\alpha\label{Ricik}
\end{eqnarray}

By (\ref{Ricik}), we immediately obtain the following

\begin{prop}
Let $F = \alpha+\beta$ be a Randers metric with ${\bf S}=0$.
Then $F$ is of zero Ricci curvature,
${\bf Ric} = 0$, 
 if and only if 
\be
s^k_{\ 0|k} = (n-1) s_j s^j_{\ 0}\ee
\be
\overline{\bf Ric}
=- (n-1) \Big (s_{0|0}+s_0s_0 \Big )
+\alpha^2 s^k_{\ j}s^j_{\ k}- 2 s_{k0}s^k_{\ 0} .
\ee
\end{prop}

\bigskip

We do not know if there is a three-dimensional Randers metric
satisfying ${\bf Ric}=0, \ {\bf S}=0$ and ${\bf K}\not=0$.
Such examples, if exist, are of interest.

\section{Examples}

In this section, we will construct two interesting examples, one in dimension two and the other in dimension three. Both metrics satisfy 
that ${\bf S} =0$ and ${\bf K}=0$. Yet, they are not locally projectively flat. More precisely, they are not Douglas metrics. 

\bigskip
\noindent{\bf Examples in Dimension two}:
Let $\alpha=\sqrt{u^2+v^2}$ denote the standard Euclidean metric on $\Bbb R^2$.
Take a vector field ${\bf v}$ on the unit disk $\Bbb D^2$ given by
\[ {\bf v} =(-y, x), \ \ \ \ \ \ p=(x, y) \in \Bbb D^2.\]
The Finsler metric associated with $(\alpha, {\bf v})$ is a Randers metric $F=\tilde{\alpha} +\tilde{\beta}$, where
\begin{eqnarray*}
 \tilde{\alpha} :& = &  {\sqrt{ (-yu+xv)^2  + (u^2+v^2) (1-x^2-y^2) } \over 1-x^2-y^2},\\
\tilde{\beta} :& = &- {-yu+xv\over 1- x^2-y^2 }.
\end{eqnarray*}

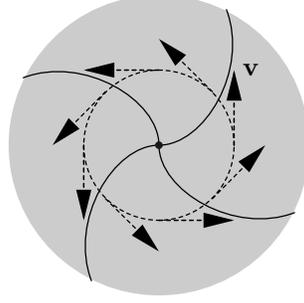
\begin{center}
\begin{texdraw}
\drawdim cm
\linewd  0.02
\lpatt(0.05 0.05)

\move(0 0)
\fcir f:0.8 r:2
\lcir r:1

\move(1 0)
\lvec(1 1)
\move(1 0)
\arrowheadtype t:F \avec(1 1)
\move(-1 0)
\lvec(-1 -1)
\move(-1 0)
\arrowheadtype t:F \avec(-1 -1)
\move(0 1)
\lvec(-1 1)
\move(0 1)
\arrowheadtype t:F \avec(-1 1)
\move(0 -1)
\lvec(1 -1)
\move(0 -1)
\arrowheadtype t:F \avec(1 -1)
\move(0.707 0.707)
\arrowheadtype t:F \avec(0 1.4142)
\move(-0.707 0.707)
\arrowheadtype t:F \avec(-1.4142 0)
\move(-0.707 -0.707)
\arrowheadtype t:F \avec(0 -1.4142)
\move(0.707 -0.707)
\arrowheadtype t:F \avec(1.4142 0)

\move(1.1 1)
\textref h:L v:C \htext{${\bf v}$}

\lpatt()
\move(0 0)
\clvec(0 0.6)(-1 1.2)(-1.8 0.9)
\move(0 0)
\clvec(-0.6 0)(-1.2 -1)(-0.9 -1.8)
\move(0 0)
\clvec(0 -0.6)(1 -1.2)(1.8 -0.9)
\move(0 0)
\clvec(0.6 0)(1.2 1)(0.9 1.8)

\move(0 0)
\fcir f:0.1 r:0.05

\move(2.2 0)
\fcir f:1 r:0.0
\move(-2.2 0)
\fcir f:1 r:0.0
\move(0 2.2)
\fcir f:1 r:0.0
\move(0 -2.2)
\fcir f:1 r:0.0
\end{texdraw}
\end{center}
\[a_{11}= {1-x^2\over \Big (1-x^2-y^2\Big )^2},
\ \ \ \ a_{12} = - {xy\over \Big ( 1-x^2-y^2\Big )^2} = a_{21},
\ \ \ \ a_{22} = { 1-y^2\over \Big (1-x^2-y^2\Big )^2 }\]
\[ b_1 = {y\over 1-x^2-y^2}, \ \ \ \ b_2 = - {x\over 1-x^2-y^2}.\]
The geodesic coefficients $\tilde{G}^1$ and $\tilde{G}^2$ of $\tilde{\alpha}$ are given by
\begin{eqnarray*}
\tilde{G}^1 & = & -{x(u^2+v^2)\over 2(1-x^2-y^2)} 
- { y(xu+yv)-v \over 1-x^2-y^2} 
\;\tilde{\beta} + {xu+yv\over 1-x^2-y^2} \; u 
\\
\tilde{G}^2 & = & - { y(u^2+v^2)\over 2(1-x^2-y^2)} 
+ {x(xu+yv)-u \over 1-x^2-y^2} \; \tilde{\beta} + {xu+yv\over 1-x^2-y^2} \;  v .
\end{eqnarray*}
The Gauss curvature $\tilde{\bf K}$ of $\tilde{\alpha}$ is given by
\[ \tilde{\bf K} =-{5+x^2+y^2\over 1-x^2-y^2}    .\]
By a direct computation, we obtain
\begin{eqnarray*}
r_{11} & = & - {2xy\over \Big (1-x^2-y^2\Big)^2}\\
r_{12} & = & {x^2-y^2\over \Big (1-x^2-y^2\Big )^2 } =r_{21}\\
r_{22} & = & {2xy\over \Big ( 1-x^2-y^2 \Big )^2 }.
\end{eqnarray*}
\begin{eqnarray*}
s_{11} & = & 0\\
s_{12} & = & {1\over \Big (1-x^2-y^2\Big )^2 } = - s_{21}\\
s_{22} & = & 0.
\end{eqnarray*}
From the above equations, we obtain the following formulas for $s_i:=b_r a^{rh}s_{hi}$.
\begin{eqnarray*}
s_1 & = & { x\over 1-x^2-y^2}\\
s_2 & = & { y\over 1-x^2-y^2}
\end{eqnarray*}
We immediately see that the following identity holds.
\[ r_{ij}= - b_is_j - b_j s_i.\]
By Proposition \ref{propS=c}, we conclude that ${\bf S}=0$.

By the above equations, we obtain the following formulas for 
$H^i := -s_j y^j y^i + \tilde{\alpha} a^{ir}s_{rl} y^l$:
\begin{eqnarray*}
H^1 & = & - { xu+yv\over 1-x^2-y^2} u - { y (xu+yv)-v \over 1-x^2-y^2} \tilde{\alpha} \\
H^2 & = & - {xu+yv\over 1-x^2-y^2} v + 
{ x (xu+yv)-u\over 1-x^2-y^2} \tilde{\alpha}.
\end{eqnarray*}
Then we obtain the following formulas for the geodesic coefficients $G^1=\tilde{G}^1+H^1$ and $G^2=\tilde{G}^2+ H^2$, 
\begin{eqnarray}
G^1 & = & - {x(u^2+v^2)\over 2( 1-x^2-y^2)}
- { y(xu+yv)- v \over 1-x^2-y^2} F \label{G1local}\\
G^2 & = & - { y (u^2+v^2) \over 2(1-x^2-y^2)} 
+ { x(xu+yv)-u \over 1-x^2-y^2} F. \label{G2local} 
\end{eqnarray}

One can easily verify that $G^1$ and $G^2$ satisfy 
\be
 S:={\pa G^1\over\pa  u} + {\pa G^2 \over\pa  v} =0.\label{GuHv}
\ee
By Proposition \ref{lemdvdv}, we know that the area form of $F$ is equal to the Euclidean area form 
\[dA_{F} =dA_{\alpha}= dx dy.\]
 By (\ref{Slocal}),  we conclude that
 ${\bf S}=0$ again.

A direct computation yields
\be
{\pa G^1 \over \pa x} + {\pa G^2\over \pa y} + {\pa G^1 \over \pa u}
{\pa G^2 \over \pa v} - {\pa G^1 \over \pa v} {\pa G^2 \over \pa u} =0.
\label{GHGH}
\ee
Plugging (\ref{GuHv}) and (\ref{GHGH}) into (\ref{RicRic}), we obtain that
${\bf Ric}=0$, hence ${\bf K}=0$.

\bigskip
One can also plug
the formulas of $G^1$ and $G^2$  in (\ref{G1local}) and (\ref{G2local}) directly into (\ref{Rik}) to verify 
that  $R^1_{\ 1}=0, R^1_{\ 2}=0, R^2_{\ 1}=0$ and $R^2_{\ 2}=0$.

\bigskip
\noindent{\bf Examples in Dimension Three or Higher}:

We first construct a Randers metric in dimension three. One can easily extend it to higher dimensions.
Let $\alpha :=\sqrt{u^2+v^2+w^2}$ denote the canonical Euclidean metric on 
$\Bbb R^3$. Take a vector ${\bf v}$ in the cylinder $\Omega:= \{ (x, y, z) \ |  \ x^2+y^2 < 1 \}$ given by
\[ {\bf v}_p : = (-y, x, 0), \ \ \ \ \  \ p=(x,y, z)\in \Omega.\]
We consider the shortest time problem under the influence of ${\bf v}$.
Image that the water in a cylindrical fish tank is rotating around the central axis, and the fishes in  the tank see the food hanging near  the water surface. Each fish has to figure out the path of  the shortest time to reach the food. The path of  shortest time 
is the shortest path of the Randers metric,  $F =\tilde{\alpha} +\tilde{\beta}$, where 
\begin{eqnarray*}
\tilde{\alpha} :& = &  { \sqrt{ (-yu+xv)^2 + (u^2+v^2+w^2)(1-x^2-y^2) }\over 1-x^2-y^2 },\\
\tilde{\beta}: & = & - {-yu+xv
\over 1-x^2-y^2 }.\\
\end{eqnarray*}

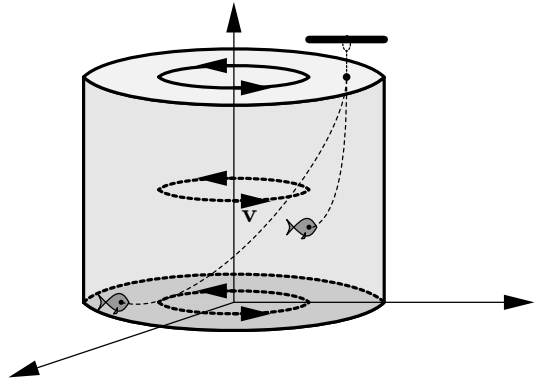
\begin{center}
\begin{texdraw}
\drawdim cm
\linewd  0.04

\move(-2 0)
\clvec(-1.5 -0.5)(1.5 -0.5)(2 0)
\lvec(2 3)
\clvec(1.5 2.5)(-1.5 2.5)(-2 3)
\lvec(-2 0)
\lfill f:0.9

\move(-2 3)
\clvec(-1.5 2.5)(1.5 2.5)(2 3)
\clvec(1.5 3.5)(-1.5 3.5)(-2 3)
\lfill f:0.95

\move(-2 0)
\lpatt(0.05 0.05)
\clvec(-1.5 0.5)(1.5 0.5)(2 0)
\lfill f:0.8
\lpatt()
\clvec(1.5 -0.5)(-1.5 -0.5)(-2 0)
\lfill f:0.8

\move(-1 3)
\clvec(-0.8 2.8)(0.8 2.8)(1 3)
\clvec(0.8 3.2)(-0.8 3.2)(-1 3) 

\move(0 2.85)
\arrowheadtype t:F \avec(0.5 2.85)
\move(0 3.15)
\arrowheadtype t:F \avec(-0.5 3.15)

\lpatt(0.05 0.05)
\move(-1 0)
\clvec(-0.8 -0.2)(0.8 -0.2)(1 0)
\clvec(0.8 0.2)(-0.8 0.2)(-1 0) 

\move(0 -0.15)
\arrowheadtype t:F \avec(0.5 -0.15)
\move(0 0.15)
\arrowheadtype t:F \avec(-0.5 0.15)

\lpatt(0.05 0.05)
\move(-1 1.5)
\clvec(-0.8 1.3)(0.8 1.3)(1 1.5)
\clvec(0.8 1.7)(-0.8 1.7)(-1 1.5) 

\move(0 1.35)
\arrowheadtype t:F \avec(0.5 1.35)
\move(0 1.65)
\arrowheadtype t:F \avec(-0.5 1.65)

\lpatt()

\linewd 0.02
\move(0 0)
\arrowheadtype t:F \avec(4 0)
\move(0 0)
\arrowheadtype t:F \avec(-3 -1)
\move(0 0)
\arrowheadtype t:F \avec(0 4)

\move(-1.4 0)
\clvec(-1.5 0.2)(-1.6 0.1)(-1.7 0)
\clvec(-1.6 -0.1)(-1.5 -0.2)(-1.4 0)
\lfill f:0.6
\move(-1.7 0)
\lvec(-1.8 0.1)
\clvec(-1.75 0)(-1.75 0)(-1.8 -0.1)
\lvec(-1.7 0)
\lfill f:0.6
\move(-1.55 -0.08)
\clvec(-1.55 -0.15)(-1.55 -0.155)(-1.6 -0.15)
\clvec(-1.6 -0.15)(-1.555 -0.1)(-1.55 -0.08)
\lfill f:0.5
\move(-1.5 0)
\fcir f:0 r:0.03

\move(1.1 1)
\clvec(1 1.2)(0.9 1.1)(0.8 1)
\clvec(0.9 0.9)(1 0.8)(1.1 1)
\lfill f:0.6
\move(0.8 1)
\lvec(0.7 1.1)
\clvec(0.75 1)(0.75 1)(0.7 0.9)
\lvec(0.8 1)
\lfill f:0.6
\move(0.95 0.92)
\clvec(0.95 0.85)(0.95 0.845)(0.9 0.85)
\clvec(0.9 0.85)(0.945 0.9)(0.95 0.92)
\lfill f:0.5
\move(1 1)
\fcir f:0 r:0.03

\move(-1.5 0)
\lpatt(0.05 0.05)
\clvec(-0.4 -0.3)(1.5 2)(1.5 3)
\move(-1.5 0)
\move(1.5 3)
\fcir f:0 r:0.05

\move(1.05 1)
\clvec(1.5 1.2)(1.5 2.5)(1.5 3)

\move(0.1 1.2)
\textref h:L v:T \htext{${\bf v}$}

\linewd 0.1
\lpatt()
\move(1 3.5)
\lvec(2 3.5)
\linewd 0.02
\lpatt(0.025 0.025)
\move(1.5 3.35)
\lvec(1.5 3)
\move(1.5 3.45)
\lellip rx:0.05 ry:0.1

\move(0 4)
\fcir f:1 r:0.0
\move(0 -1.5)
\fcir f:1 r:0.0
\end{texdraw}
\end{center}

By a similar argument, we obtain the following formulas for the geodesic coefficients $G^i$ of $F$,
\begin{eqnarray*}
G^1 & = & - {x(u^2+v^2+w^2)\over 2( 1-x^2-y^2)}
- { y(xu+yv)- v \over 1-x^2-y^2} F\\
G^2 & = & - { y (u^2+v^2+w^2) \over 2(1-x^2-y^2)} 
+ { x(xu+yv)-u \over 1-x^2-y^2} F\\
G^3 & = & 0.
\end{eqnarray*}

First, one can easily verify  that
\[ {\pa G^1 \over \pa u} + {\pa G^2\over \pa v} + {\pa G^3\over \pa w}=0.\]
By Proposition \ref{lemdvdv}, the volume form of $dV_F$ is given by
\[dV_F = dV_{\alpha} = dxdydz.\] 
By (\ref{Slocal}), we conclude that 
${\bf S}=0$.

Plugging the above formulas into (\ref{Rik}), we obtain 
$R^i_{\ k}=0$. Thus ${\bf K}=0$.

Certainly, one can extend the above Randers metric to higher dimensions with a slight modification. The formulas for  geodesic coefficients and the volume form remain same. 
Thus ${\bf S}=0$.  

Since $G^3=0$, from (\ref{Rik}), we immediately see that $R^3_{\ 1}=0, R^3_{\ 2}=0$ and $R^3_{\ 3}=0$. 
Plugging the formulas of $G^1, G^2, G^3$ into (\ref{Rik}), one can verify  that 
$R^1_{\ 1}=0, R^1_{\ 2}=0, R^1_{\ 3}=0, R^2_{\ 1}=0, R^2_{\ 2} =0$ and $R^2_{\ 3}=0$.

\noindent
Math Dept, IUPUI, 402 N. Blackford Street, Indianapolis, IN 46202-3216, USA.  \\
zshen@math.iupui.edu

\end{document}